\documentclass[11pt]{amsart}
\usepackage{amssymb,mathrsfs,bm,graphicx,geometry,mathtools,url} 
\usepackage[backref=page]{hyperref}
\usepackage[textsize=tiny]{todonotes}

\title{Statistics of centralizers of integer matrices}
\author[M. Afifurrahman]{Muhammad Afifurrahman}
\author[W. B. Wicaksono]{Wildan Bagus Wicaksono}

\address[MA]{Department of Mathematical Sciences, Ulsan National Institute of Science and Technology, 50 UNIST-gil, Ulsan 44919, South Korea}
\email{afifurrahman@unist.ac.kr}
\address[WBW]{Department of Mathematics, Universitas Brawijaya, Malang 65145, East Java, Indonesia}
\email{wildanbagus@student.ub.ac.id}

\thanks{
MA is supported by the National Research Foundation of Korea (NRF) grant funded by the Korea government (MSIT) (No. RS-2025-02262988).}

\def\cM{{\mathcal M}}

\def \C {{\mathbb C}}

\def \Q {{\mathbb Q}}
\def \R {{\mathbb R}}
\def \Z {{\mathbb Z}}

\def \bu {{\bm {u}}}
\def \bv {{\bm {v}}}

\newcommand{\rank}{\operatorname{rank}}

\newcommand{\spn}{\operatorname{span}}
\newcommand{\Cent}{\operatorname{C}}
\newcommand{\Tr}{\operatorname{Tr}}

\newtheorem{theorem}{Theorem}
\newtheorem{lemma}[theorem]{Lemma}

\newtheorem{proposition}[theorem]{Proposition}

\theoremstyle{remark}

\numberwithin{equation}{section}
\numberwithin{theorem}{section}
\numberwithin{table}{section}
\numberwithin{figure}{section}
\begin{document}

\begin{abstract} With recent interests on the statistics of commuting integer matrices, we study the related problem on the statistics of the set of centralizer of a matrix. We provide asymptotics and bounds on the size of such set when restricted to integer matrices of bounded height.
\end{abstract}
\maketitle

\section{Introduction}
\subsection{Background and results} 
This paper considers a problem on the arithmetic statistics of the set $\cM_n(\Z)$ of $n\times n$ integer matrices. Such problems are usually approached by studying the distribution of matrices in $\cM_n(\Z)$ that satisfy certain properties. 
One common way to quantify this is to study the number of matrices of bounded norm (with respect of a parameter $H$) that satisfy certain property.
Within this line, we define the set \[
\cM_n(\Z;H) \coloneqq \{A\in \cM_n(\Z) \colon |a_{i,j}|\le H \text{ for } i,j=1,\dots,n\},
\] which corresponds to the matrices in $\cM_n(\Z)$ whose sup-norm is less than $H$ and study matrices satisfying certain properties within this set, as $H\to \infty$.

There have been many recent works studying various aspects of the statistics of $\cM_n(\Z)$, such as matrices with a fixed determinant (\cite{DRS,K,Shp} for general $n\times n$ matrices,~\cite{Afif3,CM25,CM,DHP,GG1,GG2})
for $2\times 2$ matrices) or rank~\cite{K2}, or  characteristic polynomial~\cite{BGS,EMS, G2,Hab,HOS,Uri} or discriminant~\cite{OS}. Another direction is to consider the statistics of the related sets of matrices, such as the set of products of matrices~\cite{Afif1, BKOSW} and the set of commuting matrices, to which our work will mainly focus on.

This work focuses on the statistics of the set of commuting pairs of matrices. Browning--Sawin--Wang~\cite{BSW} and Chapman--Mudgal~\cite{CM} considered the statistics of the set of commuting pairs of matrices\[
\{ (A,B)\in \cM_n(\Z)^2 \colon AB=BA\};
\] see also Mudgal~\cite{Mud} for a more generalised setting.

Within this framework, we now consider the following natural extension of this result. For a matrix $A\in \cM_n(\Z)$, let the set of \textit{centralizer} of $A$, $C(A)$, defined as \[
\Cent(A) \coloneqq \{B \in \cM_n(\Q) \colon AB=BA\}.
\]  This set trivially forms a free submodule of $\cM_n(\Q)$ over $\Q$. Within our arithmetic statistics framework, we are interested in the cardinality of the set \[\Cent(A,\Z;H)\coloneqq \Cent(A)\cap \cM(\Z;H);\] that is, the number of matrices in $\Cent(A)$ whose entries are integers between $[-H,H]$. In particular, we may let $H$ be a positive integer.

Our first result provide the following formula for $\#\Cent(A,\Z;H)$  when $A\in \cM_2(\Z)$.
\begin{theorem}\label{thm:2}
    Let $A = \begin{pmatrix}
        a_1&a_2\\a_3&a_4
    \end{pmatrix}\in \cM_2(\Z)$. Also, define
    \[d:=\gcd(a_2,a_3,a_1-a_4),\quad p:=|a_1-a_4|,\quad N:=\max\{2|a_2|,2|a_3|,p\},\quad S:=\left \lfloor \frac{2dH}{N}\right \rfloor .\]
    
    Then, if $H$ is an integer,
\begin{itemize}
    \item if $A=kI_2$ for some $k\in \Z$, then $ \#\Cent(A,\Z;H)=(2H+1)^4$;
    \item else,   
    \[ \#\Cent(A,\Z;H) = (2H+1)(2S+1)-\frac{p}{d}S(S+1).\]
    In particular, as $H\to \infty$,
    \[ \#\Cent(A,\Z;H) = \frac{4d(2N-p)}{N^2}H^2 + O(H)\]
\end{itemize}
\end{theorem}
Next, we obtain the following upper bound for general $n$. 
\begin{theorem}\label{thm:n}
    Suppose that $A\in \cM_n(\Z)$, with $A\neq kI_n$ for any constant $k$. We have \[
    \#\Cent(A,\Z;H) \ll H^{n^2-2n+2}.
    \]
    Furthermore for $n\ge 4$, if $A$ is not similar to  neither $\alpha I_1 \oplus \beta I_{n-1}$ nor $\alpha I_n+\beta E_{1,2}$, we have \[
     \#\Cent(A,\Z;H) \ll H^{n^2-4n+8}.
    \]

    In addition, we also have \[
    \#\Cent(A,\Z;H)\gg H^{n},
    \] which is sharp when, for example, $A$ has $n$ distinct eigenvalues. All implied constants depend on $A$.
\end{theorem}

We are also able to classify the order of $\#\Cent(A,\Z;H)$ for the case $n=3$.
\begin{theorem}\label{thm:dim-C(A)}
    Let $A\in \mathcal{M}_3(\Z)$ be a non-scalar matrix. We have either $\#\Cent(A,\Z;H)\asymp H^3$ or $\#\Cent(A,\Z;H)\asymp H^5$. Moreover, $\#\Cent(A,\Z;H)\asymp H^5$ if and only if $A$ is similar to either $\alpha I_2\oplus \beta I_1$ $(\alpha \neq \beta)$ or $J_2(\alpha)\oplus J_1(\alpha)$, where $J_n(\alpha)$ denotes the Jordan matrix of order $n$ corresponding to to eigenvalue $\alpha$.
\end{theorem}
The main tool for Theorems~\ref{thm:n} and~\ref{thm:dim-C(A)} is Theorem~\ref{thm:subspace}, which converts the problem to counting the dimension of the related submodule $\dim \Cent(A)$. This dimension has an explicit formula related to the invariant factor of $A$, as seen in Lemma~\ref{lem:inv}. In addition, we note that the results in these theorems are also valid if we let $A\in \cM_n(\Q)$ instead of $\cM_n(\Z)$. However, this result does not need to hold for real matrices, as seen by noting that when $A= \begin{pmatrix}
    1&\sqrt{2}\\0&0
\end{pmatrix}$, \begin{equation*}
C\left( A \right) = \left\{ dI_2 \colon d\in \Q    \right\}
\end{equation*} and thus $\#C\left( A,\Z;H\right) =2H+1$.

Another interesting natural question related to these previous results is to consider the statistics related to the set of $n\times n$ matrices $A$ with given $\dim \Cent(A)$. By adapting the argument of the proof of Theorem~\ref{thm:n}, we obtain \[
n \le \dim \Cent(A)\le n^2,
\] where we consider the dimension over $\Z$ or $\Q$. Equality is attained for the lower and upper bound when (for example) $A$ has $n$ distinct eigenvalues and $A=kI$ for $k\in F$, respectively. Hence, by using known results from~\cite{HLM,OS}, most matrices $A \in \cM_n(\Z;H)$ will have $\dim \Cent(A)=n$. We will formalise this shortly.

Viewing this statement as an arithmetic statistic problem, we are naturally interested on the quantity \[
C_n(H) = \# \{A \in \cM_n(\Z;H)\colon \dim \Cent(A)>n\}.
\]
We now discuss the bounds for $C_n(H)$. For $n=2$, Theorem~\ref{thm:2} trivially implies $C_2(H)=2H+1$. 

For general $n\ge 3$, we have that $C_n(H)\le R_n(H)$, where $R_n(H)$ denotes the number of matrices in $\cM_n(\Z;H)$ with repeated eigenvalues. 
The problem of estimating $R_n(H)$ has been previously studied in the literature, see~\cite{HLM,OS}. In particular,~\cite{OS} proved \[
C_n(H)\le R_n(H) \ll H^{n^2-\Delta_n}\log H,
\] where $\Delta_n$ is a function specified in~\cite[Theorem 1.1]{OS} with $\lim_{n\to \infty} \frac{\Delta_n}{\sqrt{2n}} = 1$. 

The following result provides a lower bound on this quantity, and an asymptotic for the case $n=3$ (by using Theorem~\ref{thm:dim-C(A)}).
\begin{theorem}\label{thm:Cn}
For $n\ge 3$, we have $C_n(H)\gg H^{n+1}\log H$. Furthermore, we have $C_3(H) \asymp H^4\log H$.
\end{theorem} 



\subsection{Notations and conventions}

We use both Landau and Vinogradov notations to describe the order of an expression. For $U$ and $V$ tending to infinity with respect to some variable(s) $x$, the equivalent notations
\[
U = O(V) \quad \Longleftrightarrow \quad U \ll V \quad \Longleftrightarrow \quad V\gg U
\] 
 mean that $|U|\le c V$ for some positive constant $c$ which may depend on the dimension variables (usually $m$ and $n$), which are fixed.  Note that when $x$ represents some independent variables, we usually consider all variables tending to infinity, unless stated otherwise. Also, we sometimes put subscripts under the notation to denote dependencies.
 
 We also use the following related notations: \begin{itemize}
    \item $U=o(V)$ if and only if $\lim\limits_{U,V\to \infty} U/V \to 0$,
    \item $U\sim V$  if and only if $\lim\limits_{U,V\to \infty} U/V \to 1$,
    \item $U \asymp V$  if and only if $U\ll V$ and $V\ll U$.
\end{itemize}
For two matrices $A$ and $B$ of suitable sizes, we denote \[
A \oplus B \coloneqq \begin{pmatrix}
    A & \mathbf{0}\\ \mathbf{0} & B
\end{pmatrix},
\] with $\mathbf{0}$ denotes the zero matrix (of appropriate size). We also denote $A\sim B$ when these matrices are similar over $\C$, unless stated otherwise. All dimensions here are considered over $\Z$ or $\Q$, unless stated otherwise.

For an $n\times n$ matrix $A$, we define its $\infty$-norm as $\|A\|_\infty=\max_{1\le i,j\le n} |a_{i,j}|$.

We denote $I_n$ as the identity matrix of order $n$ and $E_{i,j}$ as the matrix (of appropriate size) whose entries are zero except of the entry on the $i$-th row and $j$-th column, which is equal to $1$. We also denote $J_n(\alpha)$ as the Jordan form matrix of order $n$ that corresponds to the eigenvalue $\alpha$.

\section{Preliminary results}
Our main tool to obtain Theorem~\ref{thm:n} is the following proposition that converts our counting problem to dimension counting problem of a submodule $S$ of $\cM_n(\Z)$ with $\dim S=k$. Note that the statement naturally holds by viewing the problem as counting integer points in the $k$-dimensional polytope defined by the basis of $S$ and computing the Ehrhart polynomial of the polytope; see~\cite[Chapter 3]{BR} for more details on this line.
However, we provide another proof here for the sake of completeness.

\begin{theorem}\label{thm:subspace}
    Let $S$ be a submodule of $\cM_n(\Z)$ with $\dim S=k$ (considered over $\Z$). Then, we have \[
    \# (S\cap \cM_n(\Z;H)) \asymp_S H^k.
    \]
\end{theorem}

\begin{proof}
    We trivially have $\# (S\cap \cM_n(\Z;H)) \ll H^k$, thus it remains to prove \[
\# (S\cap \cM_n(\Z;H)) \gg_S H^k.
    \]
The case $k=0$ is clear. Assume that $k\ge 1$, and let
$B_1,\ldots,B_k\in\mathcal{M}_n(\mathbb{Z})$ be a basis
of $S$ over $\mathbb{Z}$.
Consider matrices of the form
\[
C=d_1B_1+\cdots+d_kB_k,
\qquad d_1,\ldots,d_k\in\mathbb{Z}.
\]
Trivially, $C\in S\cap\mathcal{M}_n(\mathbb{Z})$. Then, set
\[
P:=\max_{1\le i\le k}|d_i|,
\qquad
Q:=k\max_{1\le i\le k}\|B_i\|_\infty.
\]
If $P\le\lfloor H/Q\rfloor$, then
\[
\|C\|_\infty
\le \sum_{i=1}^{k}|d_i|\|B_i\|_\infty
\le \frac{Q}{k}\sum_{i=1}^{k}|d_i|
\le PQ
\le H.
\]
Thus, $C\in S\cap\mathcal{M}_n(\mathbb{Z};H)$.

Since $B_1,\ldots,B_k$ are linearly independent, distinct
coefficient tuples $(d_1,\ldots,d_k)$ yield distinct matrices $C$.
Consequently,
\[
\#\bigl(S\cap\mathcal{M}_n(\mathbb{Z};H)\bigr)
\ge
\left(2\left\lfloor\frac{H}{Q}\right\rfloor+1\right)^k
=
\left(\frac{2H}{Q}+O(1)\right)^k
\gg_S H^k,
\]
where the implied constant depends only on the fixed basis
chosen for $S$.
\end{proof}
Note that our method remains valid if we replace $\cM_n(\Z)$ with $\cM_n(\Q)$.

We also use the following result on the dimension of $\Cent(A)$ based on its invariant factors. 

\begin{lemma}\label{lem:inv}\cite[Chapter 5, Theorem 5.15]{AW}
    Let $A\in \cM_n(F)$ be a nonzero matrix with invariant factors $q_1,\dots,q_n$, in that order, such that $q_i\mid q_{i+1}$ for $i=1,\dots,n-1$, $\deg q_i=d_i$ for $i=1,\dots,n$, and $d_1+\dots + d_n=n$ Then,  \begin{equation}\label{eqn:dimC}
    \dim \Cent(A) = \sum_{i=1}^n (2n-2i+1)d_i.
    \end{equation}
\end{lemma}
Further notions related to this lemma can be accessed at~\cite[Section~2]{Bou}.

Finally, to prove Theorem~\ref{thm:Cn}, we use the following result of Katznelson~\cite{K2} on integer matrices with given rank. Note that their argument gives asymptotical formula when the matrices are ordered in the $\ell^2$ norm, however their argument should still holds for our cases, up to the leading constant.

\begin{lemma}\label{lem:katz}\cite[Theorem 1]{K2} Let $T_n(H;k)$ denote the number of matrices in $\cM_n(\Z;H)$ of rank $k$. If $1\le k<n$, we have $T_n(H;k)\asymp H^{kn}\log H$.
\end{lemma}

\section{Proof of Theorem~\ref{thm:2}}
We easily note that if $A=kI_2$, then $\Cent(A)=\cM_2(\Q)$, which implies $ \#\Cent(A,\Z;H)=\cM(\Z;H)$.
Thus, we can now focus on the set $\Cent(A)$ with $A\ne kI$.

We first prove that each integer matrix $B\in \Cent(A)$ is uniquely represented of the form     \begin{equation}\label{eqn:Bxy}B(x,y) = \begin{pmatrix}
        y & \frac{a_2}{d}x \\ 
        \frac{a_3}{d}x & y+\frac{a_4-a_1}{d}x
    \end{pmatrix},\end{equation} where $x,y\in \Z$ and $d=\gcd(a_2,a_3,|a_1-a_4|)$.

    We first prove the uniqueness of such representation. Let $(x_1,y_1),(x_2,y_2)\in \Z^2$ such that 
    \[\begin{pmatrix}
        y_1 & \frac{a_2}{d}x_1 \\ 
        \frac{a_3}{d}x_1 & y_1+\frac{a_4-a_1}{d}x_1
    \end{pmatrix}=\begin{pmatrix}
        y_2 & \frac{a_2}{d}x_2 \\ 
        \frac{a_3}{d}x_2 & y_2+\frac{a_4-a_1}{d}x_2
    \end{pmatrix}.\]

    Comparing the upper-left entries gives $y_1=y_2$.  Since at least one of $a_2$, $a_3$, $a_1-a_4$ is nonzero, comparing the corresponding entry gives $x_1=x_2$. Hence, the representation of $B$ is unique. 

    We now prove such representation always exists. First, we recall $\Cent(A)$ forms a submodule of $\cM_n(\Z)$. Furthermore, we may see that when $A$ is not multiples of the identity matrix, $\dim \Cent(A)=2$ and $ \Cent(A)=\spn\{I_2,A\}$. Hence, for all $B\in \Cent(A)$, there exist $r,s\in \Q$ such that \begin{equation}\label{eqn:Brs}B=rI_2+sA =\begin{pmatrix}
    r+sa_1 & sa_2\\ sa_3 & r+sa_4
\end{pmatrix}.\end{equation}
Since all entries of $B$ are integers, we have 
    \[sa_2, \qquad sa_3, \qquad s(a_1-a_4)=(r+sa_1) - (r+sa_4)\]
    are all integers. Then, using Bézout's identity, there exist integers $c_1$, $c_2$, $c_3$ with 
    \[ds = c_1(sa_2) + c_2(sa_3)+ c_3[s(a_1-a_4)].\]
  
   Letting $x:=ds$ for some integer $x$ and $y:=r+sa_1$, we have $s=\frac{x}{d}$ and $r=y-\frac{a_1x}{d}$.  Substituting this to~\eqref{eqn:Brs}, we obtain~\eqref{eqn:Bxy}, which completes our proof.

   Therefore, our problem now is equivalent to counting the number of pairs of integers $(x,y)\in \Z^2$ such that $B\in \cM(\Z;H)$. Fix such $x$, we obtain $|a_2x|, |a_3x|\le dH$. Then, considering the other entries, we obtain
    \[y\in [-H,H] \cap \left [-H-\frac{a_4-a_1}{d}x, H-\frac{a_4-a_1}{d}x \right ]\]
    Since the right-hand set is non-empty  if and only if $p|x|\le 2dH$, with $p=|a_1-a_4|$, we then have
\[
N|x|=\max\{2|a_2|,2|a_3|,|a_1-a_4|\}\cdot |x| \le 2dH \iff |x| = \left \lfloor \frac{2dH}{N}\right \rfloor =S.
\]

 For each such $x$ with $|x|\le S$, the number of integers $y$ such that $B\in \cM(\Z;H)$ is exactly $2H+1-\frac{p}{d}|x|$. Therefore, we obtain
    \[ \#\Cent(A,\Z;H)=\sum_{x=-S}^S \left (2H+1-\frac{p}{d}|x|\right ) = (2H+1)(2S+1) -\frac{p}{d}S(S+1).\]
    Since $S=\frac{2dH}{N}+O(1)$, we obtain 
    \[ \#\Cent(A,\Z;H)=\frac{4d(2N-p)}{N^2}H^2+O(H).\]

\section{Proof of Theorem~\ref{thm:n}}
Based on Theorem~\ref{thm:subspace}, we may instead compute $\dim \Cent (A)$. 
We first discuss the possible values that $\dim \Cent(A)$ may take. 

Following Lemma~\ref{lem:inv}, we first note that $d_1 \le \dots \le d_n$ and $d_1+\dots +d_n=n$. Then, we easily obtain \[
n \le \dim \Cent(A) \le n^2.
\]
Equality for the lower and upper bounds are obtained if and only if $(d_1,\dots, d_n)= (0,\dots, 0,n)$ and $(1,\dots, 1)$, respectively. With respect to the matrices, we easily see that the upper bound is attainable if and only if the minimal polynomial of $A$ is linear, i.e. $A=kI$ for some constant $k$. The lower bound is attained, for example, when $A$ has $n$ different eigenvalues. This completes our lower bound proof in Theorem~\ref{thm:n}.

We now improve the upper bound for $\dim \Cent (A)$ for other cases of $A$ in the following proposition.

\begin{proposition}\label{prp:AGHM}
    If $n\ge 4$ and $A\in \cM_n(F)$ is similar to non-scalar matrix. Then,
    
    \begin{itemize}
        \item We have $\dim C(A)\le n^2-2n+2$, with equality attained if and only if $A$ is similar to either $\alpha I_1 \oplus \beta I_{n-1}$ or $\alpha I_n+\beta E_{1,2}$ for some $a,b\in F$.
        \item Otherwise, we have \[
    \dim \Cent (A)\leq n^2-4n+8.
    \]
    \end{itemize}
\end{proposition}

We note that the first part of the above proposition was first proven by  Akbari--Ghandehari--Hadian--Mohammadian~\cite[Lemma 3]{AGHM}, however here we give a different proof.
\begin{proof}
    Let $(P_1\mid \cdots\mid P_n)$ be the invariant factors of $A$ and $d_i=\deg P_i$. We have $d_1\le \cdots \le d_n$ and $\sum_{i=1}^n d_i=n$. Define 
    \[f(d_1,\cdots,d_n) := \sum_{i=1}^n (2n-2i+1)d_i.\]
    We first analyse the largest value of $f$ over all possible values of $(d_1,\dots,d_n)$.

    If all of $d_i$ are positive, we have $d_i=1$ for all $i=1,\cdots,n$. Therefore, $A$ must be a scalar matrix.

Next, we analyse the case where exactly one of $d_i$ has zero value, which implies $(d_1,d_2,\cdots,d_{n-1},d_n)=(0,1,\cdots,1,2)$ and $\dim \Cent(A)=n^2-2n+2$.

We now classify all $A$ such that this holds.
In this case, the invariant factor of $A$ must be of the form $t-\alpha,\cdots,t-\alpha,(t-\alpha)(t-\beta)$. By considering the corresponding companion matrices, we then obtain
    \[A\sim \alpha I_{n-2} \oplus \begin{pmatrix}
        0 & -\alpha\beta \\ 1 & \alpha+\beta
    \end{pmatrix}.\]
    To further classify $A$, we now find matrices that are similar to the last matrix of the above expression.
    First, if $\alpha\neq \beta$, we let $P:=\begin{pmatrix}
        -\beta & -\alpha \\ 1 & 1
    \end{pmatrix}$ and obtain $P^{-1}\begin{pmatrix}
        0 & -\alpha\beta\\ 1  & \alpha+\beta
    \end{pmatrix}P = \begin{pmatrix}
        \alpha & 0 \\ 0 & \beta
    \end{pmatrix}$. Therefore, $A\sim \alpha I_{n-1}\oplus \beta I_1$. 
    
    Next, when $\alpha=\beta$, we choose $Q:=\begin{pmatrix}
        -\alpha & 1 \\ 1 & 0
    \end{pmatrix}$ and we have $Q^{-1}\begin{pmatrix}
        0 & -\alpha^2 \\ 1 & 2\alpha
    \end{pmatrix}Q= \begin{pmatrix}
        \alpha & 1 \\ 0 & \alpha
    \end{pmatrix}=:J_2(\alpha)$. Therefore, 
    \[A\sim \alpha I_{n-2} \oplus J_2(\alpha) \sim J_2(\alpha) \oplus \alpha I_{n-2} =\alpha I_n+E_{1,2},\]
    which completes our classification.

We now consider the third-largest value of $f$. First, consider the value of $f(d_1,\dots,d_n)$ in the case where at least three of $d_i$ must be zero. Let $k$ be the index such that $d_i=0$ for all $i\le k$, and $l$ be the first index such that $d_l\ge 2$. If we transform $d_k=0$ to $d_k'=1$ and $d_l$ to $d_l'=d_l-1$, we have 
    \begin{align*}
        f(d_1,\cdots,d_k',\cdots,d_l',\cdots,d_n)-f(d_1,\cdots,d_k,\cdots,d_l,\cdots,d_n) &= (2n-2k+1) - (2n-2l+1)\\
        &= 2l-2k>0.
    \end{align*} Hence, by applying this transformation repeatedly until there are less than three of $d_i$ are zero, we obtain that $f$ obtains its third-largest value when exactly two of $d_i$ are zero. There are two such possible tuplets:
    \begin{align*}
        f(0,0,1,\cdots,1,3) = (n-2)^2+2\quad \text{and}\quad f(0,0,1,\cdots,1,2,2)=(n-2)^2+4.
    \end{align*}
    Therefore, we have that $\dim \Cent(A)\le n^2-4n+8$ for all such $A$.
\end{proof}


Our proof of Theorem~\ref{thm:n} is then complete by applying Theorem~\ref{thm:subspace} to Proposition~\ref{prp:AGHM}.








\section{Proof of Theorem \ref{thm:dim-C(A)}}
We instead prove that either $\dim \Cent(A)\in \{3,5\}$ and classify all $A\in \cM_3(\Z)$ such that $\dim \Cent(A)=5$.
    Let $P\in \cM_3(\C)$ be an invertible matrix such that $J=P^{-1}AP$ is the Jordan canonical form of $A$. We have 
    \begin{align*}
        AB = BA &\iff (P^{-1}AP)(P^{-1}BP) =(P^{-1}BP)(P^{-1}AP)\\
        &\iff J(P^{-1}BP)=(P^{-1}BP)J.
    \end{align*}
    Thereofore, $B\in C(A)$ if and only if $P^{-1}BP\in \Cent(J)$. Note that the map $f:C(A)\to \Cent(J)$ defined by $f(B)=P^{-1}BP$ is an isomorphism. Hence, $\dim C(A)=\dim \Cent(J)$, so it suffices to determine $\dim \Cent(J)$.
    
 We now divide the caseworks based on the eigenvalues of $A$ and considering a parameterised matrix $K\in \Cent(J)$, with \begin{equation}\label{eqn:JK}
     JK=KJ.
 \end{equation}
  First, suppose that $A$ has three distinct eigenvalues $\alpha$, $\beta$, $\gamma$ (that may not be real). Then, by expanding~\eqref{eqn:JK}, we have (up to block permutations)
        \[J=\begin{pmatrix}
            \alpha & 0 & 0\\ 0& \beta & 0 \\
            0& 0 & \gamma 
        \end{pmatrix}\qquad \text{and}\qquad K=\begin{pmatrix}
            k_{1,1} &0 & 0 \\ 0& k_{2,2} &0 \\ 0 & 0 &k_{3,3}
        \end{pmatrix},\]
        for parameters $k_{1,1}, k_{2,2}, k_{3,3}\in \Z$. Hence $\dim \Cent(J)=3$.
        
    Next, suppose that $A$ has exactly two distinct eigenvalues, $\alpha$ and $\beta$, with $\alpha$ having algebraic multiplicity of $2$. In this case, $A$ has two possible Jordan forms, up to block permutations.
        First, when $J=J_2(\alpha)\oplus J_1(\beta)$, by expanding~\eqref{eqn:JK}, we obtain
                \[J=\begin{pmatrix}
            \alpha & 1 & 0\\ 0& \alpha & 0 \\
            0& 0 & \beta
        \end{pmatrix}\qquad \text{and}\qquad K=\begin{pmatrix}
            k_{1,1} &k_{1,2} & 0 \\ 0& k_{1,1} &0 \\ 0 & 0 &k_{3,3}
        \end{pmatrix}.\]
        Hence, $\dim \Cent(J)=3$.

        Next, when $J=J_1(\alpha)\oplus J_1(\alpha)\oplus  J_1(\beta)$, the equation~\eqref{eqn:JK} implies
                \[J=\begin{pmatrix}
            \alpha & 0 & 0\\ 0& \alpha & 0 \\
            0& 0 & \beta
        \end{pmatrix}\qquad \text{and}\qquad K=\begin{pmatrix}
            k_{1,1} &k_{1,2} & 0 \\ k_{2,1}& k_{2,2} &0 \\ 0 & 0 &k_{3,3}
        \end{pmatrix}.\]
        Hence, $\dim \Cent(J)=5$.
        Finally, suppose that $A$ has exactly one eigenvalue $\alpha$. In this case, we have three possible Jordan forms, up to block permutations.
        First, when $J=J_3(\alpha)$,~\eqref{eqn:JK} implies
                \[J=\begin{pmatrix}
            \alpha & 1 & 0\\ 0& \alpha & 1 \\
            0& 0 & \alpha
        \end{pmatrix}\qquad \text{and}\qquad K=\begin{pmatrix}
            k_{1,1} &k_{1,2} & k_{1,3} \\ 0& k_{1,1} &k_{1,2} \\ 0 & 0 &k_{1,1}
        \end{pmatrix}.\]
        Hence, $\dim \Cent(J)=3$.

        Next, when $J=J_2(\alpha)\oplus J_1(\alpha)$,~\eqref{eqn:JK} implies
                \[J=\begin{pmatrix}
            \alpha & 1 & 0\\ 0& \alpha & 0 \\
            0& 0 & \alpha
        \end{pmatrix}\qquad \text{and}\qquad K=\begin{pmatrix}
            k_{1,1} &k_{1,2} & k_{1,3} \\ 0& k_{1,1} &0 \\ 0 & k_{3,2} &k_{3,3}
        \end{pmatrix}.\]
        Hence, $\dim \Cent(J)=5$.
        Finally, when $J=J_1(\alpha)\oplus J_1(\alpha)\oplus J_1(\alpha)=\alpha I_3$, we have  $A=PJP^{-1}=\alpha I_3$, which contradicts the hypothesis and completes our casework.
        
   Applying Theorem~\ref{thm:subspace} to each of these cases completes our proof of Theorem~\ref{thm:dim-C(A)}.

\section{Proof of Theorem~\ref{thm:Cn}}
We first prove the upper bound for $n=3$, $C_3(H)\ll H^4 \log H$. To prove this, we note that from the arguments from the proof of Theorem~\ref{thm:dim-C(A)}, any matrix $A\in \cM_3(\Z;H)$ with $\dim C(A)=5$ must have an eigenvalue $\alpha$ of multiplicity at least $2$ with $\rank(\alpha I-A)=1$. Hence, to bound $C_n(H)$, it is sufficient to count the number of matrices in $\cM_3(\Z;H)$ such that such $\alpha$ exists.

First, note that since $\alpha$ is a root of $\det(xI-\alpha)$, which is a cubic monic polynomial with integer coefficients, $\alpha$ must be an integer. Furthermore, by a classical application of Gershgorin's circle theorem, we have \[
|\alpha| \le 3H.
\] Fix one such $\alpha$, and note that $\| \alpha I-A\|_\infty \le 4H$. Then, for each possible choice of $\alpha$, we apply Lemma~\ref{lem:katz}. 
Adding over all possible $\alpha$, we have \[
C_3(H)\ll (6H+1) H^{3\cdot 1}\log H \ll H^{4}\log H,
\] which completes the proof of the upper bound of $C_3(H)$.

We now prove a lower bound of $C_n(H)$, for $n\ge 3$ with the following direct construction. Let \[A(\alpha;B)=\alpha I_n +B \in \cM_n(\Z;H),\quad \alpha \in \mathbb{Z}\] where $|\alpha|\le H/2$ and $B\in \cM_n(\Z;H/2)$ such that $\rank B=1$.

We first prove that different choices of $(\alpha;B)$ give different $A(\alpha;B)$. Let $A(\alpha_1,B_1)=A(\alpha_2,B_2)$, which implies \[
(\alpha_1-\alpha_2) I_n = B_2-B_1.
\] Note that $\rank(B_2-B_1)\le \rank B_2 +\rank (-B_1)=2$. Hence, we must have $\alpha_1=\alpha_2$, which implies $B_1=B_2$ and completes our proof.

Next, we prove that $\dim \Cent(A(\alpha;B))=n^2-2n+2$. For brevity, we write $A=A(\alpha;B)$ from now on. We first prove an identity of $B$ as follows. First, since $\rank B=1$, there exists two vectors $\bu,\bv\in \R^n$ such that $B=\bu\bv^T$ and $\Tr B = \bv^T\bu$. Let $\Tr B=t$, we then have \[
B^2= \bu\bv^T\bu \bv^T = \bu(t)\bv^T=t\bu\bv^T=tB.
\]
Furthermore, since $\rank B=1$, $B$ has an eigenvalue $0$ of algebraic multiplicity at least $n-1$. Hence, the last eigenvalue of $B$ must be $t$.

Returning to the problem of counting $\dim A$, we now divide the casework based on the value of $t$. First, when $t\ne 0$,  we have that $B\sim tE_{1,1}$. Hence, we have \[
A=\alpha I_n + B \sim \alpha I_n + tE_{1,1}  = (\alpha+t)I_1\oplus \alpha I_{n-1}.
\] Hence, by applying Proposition~\ref{prp:AGHM}, we have $\dim \Cent(A)=n^2-2n+2$  in this case.

Next, when $t=0$, we have that $B^2=0$ and all eigenvalues of $B$ are $0$. We now determine the Jordan canonical form of $B$. Since the minimal polynomial of $B$ is $x^2$, then the largest size of the Jordan block is $2\times 2$. Furthermore, because $\rank B=1$, there are exactly one  block of size $2\times 2$. Thus, the Jordan form of $B$ is  $J_{2}(0)\oplus \bm0_{n-2\times n-2} = E_{1,2}$, which implies \[ A=\alpha I_n + B \sim \alpha I_n + E_{1,2}
\] Hence, by applying Proposition~\ref{prp:AGHM}, we also have $\dim \Cent(A)=n^2-2n+2$  in this case.

Thus, in either case we have $\dim \Cent(A)=n^2-2n+2$ for any choice of $(\alpha;B)$. There are asymptotically $H$ choices for $\alpha$ and $(H/2)^{n} \log (H/2)\asymp H^n\log H$ choices for $B$, where we obtain the last fact from Lemma~\ref{lem:katz}. Hence, there are at least $H^{n+1}\log H$ matrices $A\in \cM_n(\Z;H)$ with $\dim \Cent(A)=n^2-2n+2$, which completes our lower bound proof of Theorem~\ref{thm:Cn}.


\end{document}